\documentclass[12 pt,twoside]{article}
\usepackage{graphicx}

\newtheorem{theorem}{\bf Theorem}[section]

\newtheorem{lemma}[theorem]{\bf Lemma}
\newtheorem{corollary}[theorem]{\bf Corollary}

\newtheorem{emp}[theorem]{\bf Claim}

\newtheorem{conjecture}{\bf Conjecture}

{{\footnotesize\sc }}
\input{amssym}

\begin{document}
\vspace{2 cm}
\title{{\large On three-color Ramsey number of paths}}
\vspace{4cm}
\bigskip\author{\small L. Maherani$^{\textrm{a}}$, G.R. Omidi$^{\textrm{a},\textrm{b},{1}}$, G. Raeisi$^{\textrm{a}}$, M. Shahsiah$^{\textrm{a}}$ \\
\footnotesize  $^{\textrm{a}}$ Department of Mathematical Sciences,
Isfahan
University of Technology,\\ \footnotesize Isfahan, 84156-83111, Iran\\
\footnotesize  $^{\textrm{b}}$School of Mathematics, Institute for
Research
in Fundamental Sciences (IPM),\\
\footnotesize  P.O.Box: 19395-5746, Tehran, Iran\\
\footnotesize {l.maherani@math.iut.ac.ir, romidi@cc.iut.ac.ir,
g.raeisi@math.iut.ac.ir, m.shahsiah@math.iut.ac.ir}}
\date {}
\footnotesize\maketitle\footnotetext[1] {The author was in part
supported by a grant from IPM (No. 90050049)} \footnotesize

\begin{abstract}\rm{}
\medskip
\footnotesize Let $G_1, G_2, \ldots , G_t$ be graphs. The
multicolor Ramsey number  $R(G_1, G_2, \ldots, G_t)$ is the
smallest positive integer $n$ such that if the edges of complete
graph $K_n$ are partitioned into $t$ disjoint color classes giving
$t$ graphs $H_1,H_2,\ldots,H_t$, then at least one $H_i$  has a
subgraph isomorphic to $G_i$. In this paper, we prove that if
$(n,m)\neq (3,3), (3,4)$ and $m\geq n$, then
$R(P_3,P_n,P_m)=R(P_n,P_m)=m+\lfloor \frac{n}{2}\rfloor-1$.
Consequently $R(P_3,mK_2,nK_2)=2m+n-1$ for $m\geq n\geq 3$.
\vspace{.5cm}
 \\{ {Keywords}:{ \footnotesize Ramsey Number, Path.\medskip}}
\noindent
\\{\footnotesize {AMS Subject Classification}:  05C55.}
\end{abstract}
\small
\medskip
\section{\normalsize{Introduction}}

\medskip
\hspace{.5 cm}In this paper, we only concerned with undirected simple finite graphs and we follow \cite{Boundy} for terminology
and notations not defined here. For a graph $G$, the vertex set, edge set, maximum degree and
minimum degree of $G$ are denoted by $V(G)$, $E(G)$, $\Delta(G)$ and
$\delta(G)$ (or simply $V$, $E$, $\Delta$, $\delta$), respectively. As usual, the complete graph of order $p$
 is denoted by $K_p$ and a complete bipartite graph with partite set $(X,Y)$ such that $|X|=m$ and $|Y|=n$ is
  denoted by $K_{m,n}$. For two disjoint subsets $X$ and $Y$ of the vertices of a graph $G$, we use $E(X,Y)$ to denote the set of
  all edges with one end point in $X$ and the other in $Y$. For a vertex $v$ and an induced subgraph $H$ of $G$  the set of all neighbors of
  $v$ in $H$ are denoted by $N_H(v)$.  Throughout this paper, we denote a cycle and a path on $m$ vertices
   by $C_m$ and $P_m$, respectively.
   Also for a 3-edge coloring (say green, red and blue) of a graph $G$, we denote by $G^g$ (resp. $G^r$ and $G^b$) the induced subgraph
    by the edges of color green (resp. red and blue).

\bigskip For given graphs $G_1, G_2, \ldots , G_t$ the {\it multicolor Ramsey number}
$R(G_1, G_2, \ldots, G_t)$, is the smallest positive integer $n$
such that if the edges of  complete graph $K_n$ are partitioned
into $t$ disjoint color classes giving $t$ graphs
$H_1,H_2,\ldots,H_t$, then at least one $H_i$  has a subgraph
isomorphic to $G_i$.
 The existence of such a positive integer is guaranteed by Ramsey's classical result \cite{Ramsey}.
  Since 1970's, Ramsey theory has grown into one of the most active
   areas of research within combinatorics, overlapping variously with graph theory, number theory, geometry and logic.
For $t\geq 3$, there is a few results about multicolor Ramsey number $R(G_1, G_2, \ldots, G_t)$.
 A survey including some results on Ramsey number of graphs, can be found in \cite{survey}.

\bigskip The multicolor Ramsey number $R(P_{n_1}, P_{n_2},\ldots, P_{n_t})$ is not known for $t\geq 3$.
 In the case $t=2$, a well-known theorem of Gerencs\'{e}r and Gy\'{a}rf\'{a}s \cite{Ramsey number of paths} states that
  $R(P_n,P_m)=m+\Big\lfloor\frac{n}{2}\Big\rfloor-1,$
where $m\geq n\geq 2$. Faudree and Schelp in \cite{Path Ramsey
numbers in multicolorings} determined $R(P_{n_1},
P_{2n_2+\delta},\ldots, P_{2n_t})$ where $\delta\in \{0,1\}$ and
$n_1$ is sufficiently large. As an  improvement of this result
in \cite{Omidi Raeisi} the authors determined $R(C_{n_1},
P_{2n_2+\delta},\ldots, P_{2n_t})$ where $\delta\in \{0,1\}$ and
$n_1$ is sufficiently large. In addition, in  \cite{Path Ramsey
numbers in multicolorings} the authors determined $R(P_{n_1},P_{n_2},P_{n_3})$
for the case $n_1\geq 6(n_2+ n_3)^2$ and they conjectured that
\medskip

$$
R(P_n,P_n,P_n)= \left\lbrace
\begin{array}{ll}
2n-1  &\mbox{if~}  n~\mbox{is~odd},\vspace{.5 cm}\\
2n-2  &\mbox{if~}  n~\mbox{is~even}.
\end{array}
\right.\vspace{.2 cm}
$$

\medskip
\noindent This conjecture was established by Gy\'{a}rf\'{a}s et al. \cite{Ramsey number of three-paths} for sufficiently large $n$.
 In asymptotic form, this was proved by Figaj and
Luczak in \cite{Ramsey number for a triple of long even cycles} as
a corollary of more general results about the asymptotic results
on the Ramsey number for three long even cycles.

\bigskip  It is a natural question to ask whether similar conclusion is true
if $K_{R(P_m,P_n)}$ is replaced by some weaker structures. One
such result was obtained in \cite{tripartite Ramsey number of
paths} where it was proved that in every 2-coloring of the edges
of the complete 3-partite graph $K_{n,n,n}$ there is a
monochromatic $P_{(1-o(1))2n}$. The following conjecture involving
the minimum degree, was formulated by Schelp \cite{schelp conj}.

\medskip
\begin{conjecture}
Suppose that $n$ is large enough and $G$ is a graph on
$R(P_n,P_n)$ vertices with minimum degree larger than
$\frac{3}{4}|V(G)|$. Then in any 2-coloring of the edges of $G$
there is a monochromatic $P_n$.
\end{conjecture}

\medskip
Schelp also noticed that the condition on the minimum degree is
sharp. Indeed, suppose that $3n-1=4m$ and consider a graph whose
vertex set is partitioned into four parts $A_1,A_2,A_3,A_4$ with
$|A_i|=m$. There are no edges from $A_1$ to $A_2$ and from $A_3$
to $A_4$. Edges between $A_1,A_3$ and $A_2,A_4$ are red, edges
between $A_1,A_4$ and $A_2,A_3$ are blue and for $i=1,2,3,4$ the
edges of with two end points in $A_i$ are colored arbitrary. In
this coloring the longest monochromatic path has $2m$ vertices,
much smaller then $2n$, while the minimum degree is
$\frac{3}{4}|V(G)|-1$. Thus, this makes the conjecture surprising,
even a minuscule increase in the minimum degree results in a
dramatic increase in the length of the longest monochromatic path.
Schelp \cite{schelp} proved that there exists a $c<1$ for which
Conjecture 1 holds if the minimum degree is raised to $c|V(G)|$.
 The main result of this paper is the following.
\begin{theorem}\label{main result}
If $m\geq n$ and $(n,m)\neq(3,3),(3,4)$, then
$R(P_3,P_n,P_m)=m+\lfloor \frac{n}{2}\rfloor-1$. Moreover,
$R(P_3,P_3,P_3)=R(P_3,P_3,P_4)=5$.
\end{theorem}

In other words, $R(P_3,P_n,P_m)=R(P_n,P_m)$ for $m\geq n$ and
$(n,m)\neq (3,3), (3,4)$. Clearly $R(P_n,P_m)$ is a lower bound
for $R(P_3,P_n,P_m)$ and so we shall always prove just the claimed
upper bound for the Ramsey number.

\medskip
\section{\normalsize {$R(P_3,P_n,P_m)$ for $m\geq n$ and $n\leq 7$ }}

\medskip
\hspace{.5 cm}In this section, we provide the exact values of $R(P_3,P_n,P_m)$ when $3\leq n\leq7$ and $m\geq n$. First, we recall a result of Faudree and
Schelp.

\medskip
\begin{theorem}{\rm(\cite{Path Ramsey numbers in multicolorings})}\label{faudree and schelp}
If $G$ is a graph with $|V(G)|=nt+r$ where $0\leq r < n$ and $G$ contains no path on $n+1$ vertices, then $|E(G)|\leq t{n\choose 2}+{r\choose 2}$
with equality if and only if  either $G\cong tK_n\cup K_r$ or if $n$ is odd, $t>0$ and
 $r=(n\pm 1)/2$  $$G\cong lK_n\cup \Big(K_{{(n-1)}/{2}}+\overline{K}_{{((n+1)}/{2}+(t-l-1)n+r)}\Big),$$ for some $0\leq l<t$.
\end{theorem}

\medskip
By Theorem \ref{faudree and schelp}, it is easy to obtain the following corollary.

\medskip
\begin{corollary}\label{extermal P5,P4}
For all integer $n\geq 3$,\\
\begin{eqnarray*}
&&ex(n,P_4)=\left\lbrace
\begin{array}{ll}
n  &\mbox{if~}  n=0~\mbox{(mod~3)},\vspace{.2 cm}\\
n-1  &\mbox{if~}  n=1,2~\mbox{(mod~3)}.
\end{array}
\right.
\\
\\ \vspace{1.9 cm}
&&ex(n,P_5)=\left\lbrace
\begin{array}{ll}
3n/2  &\mbox{if~}  n=0~\mbox{(mod~4)},\vspace{.2 cm}\\
3n/2-2  &\mbox{if~}  n=2~\mbox{(mod~4)},\vspace{.2 cm}\\
(3n-3)/2 &\mbox{if~}  n=1, 3~\mbox{mod~4}.
\end{array}
\right.
\\
\\
&&ex(n,P_6)=\left\lbrace
\begin{array}{ll}
2n  &\mbox{if~}  n=0~\mbox{(mod~5)},\vspace{.2 cm}\\
2n-2  &\mbox{if~}  n=1, 4~\mbox{(mod~5)},\vspace{.2 cm}\\
2n-3 &\mbox{if~}  n=2, 3~\mbox{mod~5}.
\end{array}
\right.
\end{eqnarray*}
\end{corollary}

\medskip
\begin{theorem}{\rm(\cite{Electron. J. Combin,Discuss. Math. Graph
Theory})}\label{R(P3,P4,Pm) and R(P3,P5,Pm)}
 $R(P_3,P_4,P_m)=m+1$ for $m\geq 6$ and
$R(P_3,P_5,P_m)=m+1$ for $m\geq 8$.
\end{theorem}

\begin{theorem}\label{R(P3,P3,Pm)}
{\rm(i)}~$R(P_3,P_3,P_m)=m$ for $m\geq 5$. \medskip

\noindent{\rm(ii)}~$R(P_3,P_4,P_m)=m+1$ for $4\leq m\leq 5$.\medskip

\noindent{\rm(iii)}~$R(P_3,P_5,P_m)=m+1$ for $5\leq m\leq 7$.
\end{theorem}
\noindent\textbf{Proof. }(i) Let $G=K_m$ be 3-edge colored green, red and blue such that $G$ does not
 contain green or red $P_3$. It is clear to see that $G^b$ is connected and $\delta(G^b)\geq m-3$. Thus $G^b$ has
  a Hamiltonian path(see \cite{Boundy}) and so  a $P_m$ .
\bigskip

\noindent (ii) Let $G=K_{m+1}$ be 3-edge colored green, red and
blue such that $P_3\nsubseteq G^g$ and $P_4\nsubseteq G^r$. First
let $m=4$. Using corollary \ref{extermal P5,P4} we may assume that
$|E(G^g)|\leq 2$ and $|E(G^r)|\leq 4$. If $|E(G^r)|=4$, then by
Theorem \ref{faudree and schelp} $G^r\cong K_3\cup K_2$ or
$G^r\cong K_{1,4}$ which clearly the complement of $G^r$ with
respect to $G$ is colored green and blue and so it contains  a
blue copy of $P_4$. Thus we may assume that $|E(G^r)|\leq 3$ and
so $|E(G^b)|\geq 5$. Using corollary \ref{extermal P5,P4} $G^b$
contains $P_4$. By a similar argument one can show that
$R(P_3,P_4,P_5)=6$.
\bigskip

\noindent (iii) Let $G=K_{m+1}$ be 3-edge colored green, red and blue such that $P_3\nsubseteq G^g$
 and $P_5\nsubseteq G^r$.
First let $m\neq 5$. By a result in \cite{survey},
$R(P_3,C_4,P_m)=m+1$ for $m\in \{6,7\}$ and so we may assume that
$G$ contains a red $C_4$. Set $A=V(C_4)$ and $B=V(G)\setminus A$.
Since $P_5\nsubseteq G^r$, all edges between $A$ and $B$ are
colored green or blue which clearly $G[E(A,B)]$ contains a blue
$P_m$. Now consider the case $m=5$. By a similar
 argument, we may assume that $G^r$ and $G^b$ don't contain $C_4$ as subgraph. Since $|E(G)|=15$, by Theorem
 \ref{faudree and schelp} we may assume that $|E(G^g)|=3$, $|E(G^r)|=6$ and $|E(G^b)|=6$ and so the green edges
 form a perfect matching. But $R(P_3,P_4,P_5)=6$, by part (ii), and so we may assume that $G^r$ contains a copy
 of $P_4$, say $P=v_1v_2v_3v_4$. Set $A=V(G)\setminus V(P)=\{v_5,v_6\}$. Since $P_5\nsubseteq G^r$,
 all edges in $E(\{v_1,v_4\},A)$ are colored green or blue. Also since the green edges form a perfect matching,
  the subgraph of $G^g$ induced by $E(\{v_1,v_4\},A)$ dose not contain a perfect matching. Thus we may
   assume that $P^{\prime}=v_5v_1v_6v_4\subseteq G^b$ and $v_4v_5\in E(G^g)$. Now since $P_5\nsubseteq G^r$, at
   least one of $v_2v_5$ or $v_3v_5$, say $v_2v_5$, must be blue and so $v_3v_5P'v_4$ form a blue
   $P_5$. This observation completes the proof.
$\hfill\blacksquare$
\\

Combining Theorems \ref{R(P3,P4,Pm) and R(P3,P5,Pm)} and
\ref{R(P3,P3,Pm)}, we obtain that $R(P_3,P_n,P_m)=R(P_n,P_m)$ if
$m\geq n$, $n\in\{3,4,5\}$ and $(n,m)\neq (3,3),(3,4)$. In the
rest of this section we prove that $R(P_3,P_n,P_m)=R(P_n,P_m)$ for
$m\geq n$, $n\in\{6,7\}$. But before that we need some lemmas.
\medskip
\begin{lemma}\label{k3,4}
Let $G$ be a graph obtained from the complete bipartite graph
$K_{3,4}$ by removing an edge.
 If each edge of $G$ is colored red or blue, then $G^r$ contains $P_3$ or $G^b$ contains $P_7$.
\end{lemma}
\noindent\textbf{Proof. }Let $G=(X,Y)$, $X=\{x_1,x_2,x_3\}$ and $Y=\{y_1,y_2,y_3,y_4\}$. Also let $x_1y_1$ be the
edge that removed from $K_{3,4}$. If $G^r$ does not contain $P_3$, then $G^r$ has at most three edges. Let $H$ be a spanning subgraph of $G$
with  $E(H)=E(G^r)\cup\{x_1y_1\}$.
 It is clear to see that $H\subseteq P_3\cup 2P_2$ or $H\subseteq P_4\cup P_2\cup P_1$
and so the complement of $H$ with respect to $K_{3,4}$ contains a
copy of $P_7$.
 This observation completes the proof.
$\hfill \blacksquare$



\begin{lemma}\rm{}\label{main lemma of R(P3,P7,Pm)}
Suppose $m\geq 7$ and  the edges of $K_{m+2}$ are colored with
colors green, red and blue such that $G^b$ contains a copy of
$P_{m-1}$ as a subgraph. Then $K_{m+2}$ contains a green $P_3$, a
red $P_7$ or a blue $P_m$.
\end{lemma}
\textbf{Proof. }Assume that $G=K_{m+2}$ with
$V(G)=\{v_1,v_2,\ldots, v_{m+2}\}$ and $P=v_1v_2\ldots v_{m-1}$ is
the desired copy of $P_{m-1}$ in $G^b$. We suppose that $G^b$
contains no copy of $P_m$, then we prove that $K_{m+2}$ contains
a green $P_3$ or a red $P_7$.
We find two vertices $v,v^{\prime}\in P_{m-1}$ such that the
bipartite graph with parties $X=\{v_{m},v_{m+1},v_{m+2}\}$ and
$Y=\{v_1,v,v^{\prime},v_{m-1}\}$ is a red-green graph with at
least 11 edges and then we use Lemma \ref{k3,4}, which guarantees
the existence of a green $P_3$ or a red $P_7$. Note that we may
assume that in $G^b$, the vertices $v_2$ and $v_{m-2}$ don't have
a common neighbor in $X$. Otherwise, since $P_m\nsubseteq G^b$,
$v_3$ (also $v_{m-3}$) is not adjacent to any vertex of $X$ in
$G^b$ and so $v_3$ and $v_{m-3}$ are the desired vertices. Thus we
may assume that in $G^b$ one of $v_2$ or $v_{m-2}$, say $v_{m-2}$,
has at most one neighbor in $X$. If $N_{G^b}(v_{m-2})\cap
X=\emptyset$, then we may assume that in $G^b$ each vertex $v_i\in
V(P)\setminus \{v_1,v_{m-2},v_{m-1}\}$ has at least two neighbors
in $X$, otherwise set $v=v_i$ and $v^{\prime}=v_{m-2}$. Therefore
if $N_{G^b}(v_{m-2})\cap X=\emptyset$ we have $N_{G^b}(v_2)\cap
N_{G^b}(v_{3})\cap X \neq \emptyset$, and so $P_m\subseteq G^b$, a
contradiction. Hence $|N_{G^b}(v_{m-2})\cap X|=1$,

\bigskip
Since $|N_{G^b}(v_{m-2})\cap X|=1$, so we may assume that in $G^b$
each vertex $v_i\in V(P)\setminus \{v_1,v_{m-2},v_{m-1}\}$ has at
least one neighbor in $X$, otherwise set $v=v_i$ and
$v^{\prime}=v_{m-2}$. Since $P_m\nsubseteq G^b$, one can easily
check that $|N_{G^b}(v_i)\cap X|=1$, $2\leq i\leq 5$ and w.l.g
$N_{G^b}(v_{2})\cap X=\{v_m\}$, $N_{G^b}(v_{3})\cap
X=\{v_{m+1}\}$, $N_{G^b}(v_{4})\cap X=\{v_{m+2}\}$ and
$N_{G^b}(v_{5})\cap X=\{v_{m+1}\}$. If $m=7$ then
$v_1v_2v_3v_8v_5v_4v_9$ is a blue $P_7$ in $K_9$, a contradiction.
Now let $m\geq 8$. Since $|N_{G^b}(v_{m-2})\cap X|=1$, $v_{m-2}$
must be adjacent to a vertex in $X$ by blue and in any case we
have a copy of $P_m\subseteq G^b$, a contradiction. This
observation completes the proof. $\hfill\blacksquare$

\medskip
\begin{lemma}\rm{}\label{R(P3,P6,P7)}
$R(P_3,P_6,P_7)=9$.
\end{lemma}
\textbf{Proof. }Let $G=K_9$ be 3-edge colored  with colors green,
red and blue. By a result in \cite{Euro. J. Combin},
$R(P_3,C_6,C_6)=9$ and so we may assume that $G^b$ contains a copy
of $C_6$ as subgraph. Set $X=V(K_9)\setminus V(C_6)$.
 We may assume that all edges between $X$ and $C_6$ are colored red or green. Therefor by Lemma \ref{k3,4}, $K_9$ must contain
   a green $P_3$ or a red $P_6$, which completes the proof.
$\hfill\blacksquare$

\bigskip
Using Lemmas \ref{main lemma of R(P3,P7,Pm)} and \ref{R(P3,P6,P7)}
we have the following.

\medskip
\begin{theorem}\label{R(p3,Pm,P6),R(P3,Pm,P7)}
$R(P_3,P_7,P_m)=R(P_3,P_6,P_m)=m+2$ for $m\geq 7$. Moreover
$R(P_3,P_6,P_6)=8$.
\end{theorem}
\noindent\textbf{Proof. }Since $R(P_3,P_6,P_m)\leq R(P_3,P_7,P_m)$
and $m+2=R(P_6,P_m)\leq R(P_3,P_6,P_m)$, it is sufficient to show
that $R(P_3,P_7,P_m)\leq m+2$ for $m\geq 7$. Using Lemmas
\ref{main lemma of R(P3,P7,Pm)} and \ref{R(P3,P6,P7)} we have
$R(P_3,P_7,P_7)=9$ and again using Lemma \ref{main lemma of
R(P3,P7,Pm)} and induction on $m$ we obtain that
$R(P_3,P_7,P_m)\leq m+2$. On the other hand $8=R(P_6,P_6)\leq
R(P_3,P_6,P_6)$. To complete the proof it is sufficient to show
that $R(P_3,P_6,P_6)\leq 8$. Let $G=K_8$ be 3-edge colored with
colors green, red and blue. Suppose $G$ have neither a green $P_3$
nor a blue $P_6$. If $G$ has a red $P_6$ we are done. So suppose
that $G$ does not have any red $P_6$. Using $(iii)$ of Theorem
\ref{R(P3,P3,Pm)} we may assume that $G$ has a red $P_5$ with
vertices $v_1, v_2,\cdots, v_5$ as a subgraph. Then we may assume
that $v_1v_6, v_1v_7, v_5v_7, v_5v_8$ are blue edges. If
$E(\{v_6,v_8\},\{v_2,v_3,v_4\})$ has a blue edge, combining this
edge with the path $v_6v_1v_7v_5v_8$ gives a blue $P_6$, a
contradiction. So $v_8v_2, v_8v_4, v_6v_2, v_6v_4$ are red edges
and $v_3v_8$ and $v_3v_6$ are green edges and hence $G$ has a
green $P_3$, a contradiction.


$\hfill\blacksquare$

\medskip
\section{\normalsize{$R(P_3,P_n,P_m)$ for $m\geq n \geq 8$}}
\medskip

\hspace{.5 cm}In this section, we compute the value of
$R(P_3,P_n,P_m)$ for $m\geq n \geq8$. Before that we need some
lemmas.

\medskip
\begin{lemma}\label{red path}
Suppose that $G=K_{m+\lfloor \frac{n}{2}\rfloor-1}$, $m\geq n
\geq8$, is 3-edge colored green, red and blue and $P=v_1v_2\cdots
v_{m-1}$ is the maximum path in $G^b$. Let $A=V(G)\backslash V(P)$
and $H$ be the subgraph of $G^r$ induced by the edges in
$E(V(P)\backslash \{v_1,v_{m-1}\},A)$. Then either $P_3\subseteq
G^g$ or $d_{H}(v_i)\geq2$ for some $i$, $2\leq i\leq m-2$.
\end{lemma}
\textbf{Proof. }We suppose that $G^g$ contains no copy of $P_3$.
Since $m\geq n \geq8$, we obtain that $|A|\geq 4$ and so let
$X=\{u_1,u_2,u_3,u_4\}\subseteq A$. Again since $m\geq 8$ there is
a $v_j\in V(P)\backslash \{v_1,v_{m-1}\}$ such that all edges in
$E(X,\{v_j\})$ are red and blue.
If $|N_{G^r}(v_j)\cap X|\geq 2$, we have nothing to prove.
Otherwise, we may assume that $Y=\{u_1,u_2,u_3\}\subseteq N_{G^b}(v_j)\cap
X$. Since $P_m\nsubseteq G^b$, $G^r$ contains at
least two edges in  $E(Y,\{v\})$ for some $v\in
\{v_{j-1},v_{j+1}\}$
 and so $d_H(v)\geq 2$.
$\hfill\blacksquare$

\medskip
\begin{lemma}\label{independent red edges}
Suppose that $G=K_n$ is 3-edge colored green, red and blue,
$P_3\nsubseteq G^g$ and $P$ is a maximal path in $G^b$ with
endpoints $x$ and $y$. Then for every two vertices $z$ and $w$ of
$V(G)\setminus V(P)$ either $xz,yw\in E(G^r)$ or $xw,yz\in
E(G^r)$.
\end{lemma}
\textbf{Proof. }Since $P_3\nsubseteq G^g$ and $P$ is a maximal
path in $G^b$, each of  $z$ and $w$ is adjacent  to at least one
of $x$ and $y$ in $G^r$. With no loss of generality, suppose that
$xz\in E(G^r)$. If $yw\in E(G^r)$, the proof is completed.
Otherwise, $yw\in E(G^g)$ and so $xw$, $yz \in E(G^r)$, which
completes the proof.

 $\hfill\blacksquare$

\begin{lemma}\label{P3,P8,P8}
$R(P_3,P_8,P_8)=11$
\end{lemma}
\textbf{Proof. }Let $G=K_{11}$ be 3-edge colored green, red and
blue such that $P_3\nsubseteq G^g$. We find monochromatic copy of
$P_8$ in blue or red color. By Theorem
\ref{R(p3,Pm,P6),R(P3,Pm,P7)}, $R(P_3,P_7,P_8)=10$ and so we may
assume that $P_7$ is a maximum path in $G^r$.  Let $P=v_1v_2\ldots
v_7\subseteq G^r$ and $A=V(G)\backslash V(P)=\{x_1,x_2,x_3,x_4\}$.
Using Lemma \ref{red path}, there exists a $v_j\in V(P)\backslash
\{v_1,v_7\}$ which is adjacent to at least two vertices of $A$,
say $x_1, x_2$, in $G^b$. By Lemma \ref{independent red edges},
w.l.g  we may assume that
$\{x_1v_1,x_2v_7,x_3v_1,x_4v_7\}\subseteq E(G^b)$ and so
$Q_7=x_3v_1x_1v_jx_2v_7x_4\subseteq G^b$. Let $K=V(P)\backslash
\{v_1,v_j,v_7\}$. Then $|K|=4$ and one can easily check that at
least one of $x_3$ or $x_4$ is adjacent to a vertex of $K$, say
$v_i$, in $G^b$. Therefore $Q_7\cup \{v_i\}$ is a blue $P_8$.
$\hfill\blacksquare$

\hspace{1.6 cm}
\begin{theorem}\label{R(P3,Pm,Pn)}
For any  $m\geq n \geq8$, $R(P_3,P_n,P_m)=m+\lfloor
\frac{n}{2}\rfloor-1$.
\end{theorem}
 \textbf{Proof. } Let $t=m+\lfloor \frac{n}{2}\rfloor-1$ and $G=K_t$ be 3-edge colored green,
  red and blue such that $P_3\nsubseteq G^g$ and  $P_m\nsubseteq G^b$. By induction on $m+n$, we
  prove that  $P_n\subseteq G^r$. By Lemma
 \ref{P3,P8,P8} theorem is true for $m=n=8$. By the induction hypothesis $R(P_3,P_n,P_{m-1})\leq m+\lfloor
\frac{n}{2}\rfloor-1$
 and so there is a
$P_{m-1}\subseteq G^b$. Let $P=P_{m-1}=v_1v_2\ldots v_{m-1}$,
$A=V(G)\backslash V(P)$ and $H$ be the subgraph  of $G^r$ induced
by the edges in $E(V(P)\backslash \{v_1,v_{m-1}\},A)$. Suppose $Q$
is a maximal path of $H$ with end points $u_1$  and $u_2$ in $A$,
the existence of such a path is guaranteed by Lemma \ref{red
path}.
 Let $K=(V(P)\setminus \{v_1,v_{m-1}\})\setminus V(Q)$. If all
vertices in $A$ are covered by $Q$, then by Lemma \ref{independent
red edges}, we may assume that $u_1v_1, u_2v_{m-1}\in E(G^r)$ and
so $R=v_1u_1Qu_2v_{m-1}$ is a red path on
$2\lfloor\frac{n}{2}\rfloor+1$ vertices. Thus we may assume that
$A\setminus V(Q)\neq \emptyset$.

\vspace{.6 cm} \noindent\textbf{Case 1. }$|A\backslash
V(Q)|=1$.\vspace{.3 cm}

\noindent Let $A\backslash V(Q)=\{x\}$. By Lemma \ref{independent
red edges}, we may assume that $v_1u_1, v_{m-1}u_2\in E(G^r)$. In
the other hand, since $P$ is maximal and $P_3\nsubseteq G^g$, $x$
is adjacent to at least one of $v_1$ and $v_{m-1}$ in $G^r$, say
$v_1$. Thus $R=xv_1u_1Qu_2v_{m-1} \subseteq G^r$ form a  path on
$2\lfloor\frac{n}{2}\rfloor$ vertices. If $n$ is even, there is
nothing to prove and so we may assume that $n$ is odd. Note that
$|K|=m-3-(\lfloor\frac{n}{2}\rfloor-2)\geq \lceil
\frac{m}{2}\rceil-1>\lfloor\frac{n}{2}\rfloor-1$ and so by the
Pigeonhole principle there exist two consecutive vertices
$v_i,v_{i+1}$ in $K$. If $xv_i \in E(G^r)$ (or $xv_{i+1}\in
E(G^r))$, then $\{v_i\}\cup V(R)$ (or $\{v_{i+1}\}\cup V(R))$ form
a red $P_n$. Otherwise, since both $xv_i$ and $xv_{i+1}$ are not
in $E(G^g)$ or $E(G^b)$, w.l.g we may assume that $xv_i\in E(G^b)$
and $xv_{i+1}\in E(G^g)$ which implies that $xv_{m-1}\in E(G^r)$.
Therefore $V(R)\cup \{x\}$ form a copy of $C_{n-1}$ in $G^r$. It
is clear to see that at least one of $v_i$ or $v_{i+1}$ is
adjacent  to one of $u_1$ or $u_2$ by a red edge Thus, we can find
a red $P_n$.

\vspace{.6 cm} \noindent\textbf{Case 2. }$|A\setminus
V(Q)|=2$.\vspace{.3 cm}

\noindent Let $A\setminus V(Q)=\{x,y\}$. Using Lemma
\ref{independent red edges} we may assume that
$R=xv_1u_1Qu_2v_{m-1}y$ is a red path
 on $2\lfloor\frac{n}{2}\rfloor-1$ vertices. $($Note that
in this case, $|K|=m-3-(\lfloor\frac{n}{2}\rfloor-3)\geq \lceil
\frac{m}{2}\rceil\geq\lceil\frac{n}{2}\rceil)$. We consider the
 following subcases.\\

\noindent {\it Subcase 1}. $n$ is even:

\medskip \noindent By the Pigeonhole principle there exists a
pair of vertices $(v_i,v_{i+1})$ in $K$. If one of $x$ or $y$ is
adjacent to one of $v_i$ or $v_{i+1}$, say $v_i$, in $G^r$, then
$v_ixRy$ form a red $P_n$. Otherwise,  green and also
blue edges in $E(\{v_i,v_{i+1}\},\{x,y\})$ form a matching and so
$yv_1$ is red and w.l.g we may assume that $u_1v_i$ is red. Thus
$R'=v_iu_1Qu_2v_{m-1}yv_1x\subseteq G^r$ is a path on $n$ vertices.\\

\noindent {\it Subcase 2}. $n$ is odd:

\medskip \noindent By
the Pigeonhole principle there exist two disjoint pairs of
vertices $(v_j,v_{j+1})$ and $(v_k,v_{k+1})$ in $K$. It is easy to
see that each of $x$ and $y$ is adjacent to a vertex in
$B=\{v_j,v_{j+1},v_k,v_{k+1}\}$ by red edge. If the mentioned
neighbors of $x$ and $y$ are distinct we have a red $P_n$,
otherwise let $v_j\in B$ be the only neighbor of $x$ and $y$.
Therefore, $\{v_j\}\cup V(R)$ form a red $C_{n-1}$. It
is easy to see that there is an edge in $G^r$ between $B\backslash
\{v_j\}$ and $\{u_1,u_2\}$ and so a red $P_n$ can be found.

\vspace{.6 cm} \noindent\textbf{Case 3. }$|A\setminus V(Q)|\geq
3$.\vspace{.3 cm}

\noindent Let $x,y,z\in A\setminus  V(Q)$.
\begin{emp}\label{red path Q'}
Let $H$ be the subgraph of $G^r$ induced by the edges in
$E(A\backslash V(Q),K)$. There is a vertex $v\in H\cap K$ such
that $d_H(v)\geq 2$.

\end{emp}
{\bf Proof.}
 There are at least
$\lceil\frac{n}{2}\rceil+1$ vertices in $K$. By the Pigeonhole
principle, there are two disjoint pairs of vertices
$(v_i,v_{i+1})$ and $(v_j,v_{j+1})$ in $K$. We prove the claim by
considering the number of red edges from $\{v_i,v_{i+1}\}$ to
$\{x, y,z\}$. If there are more than two such edges, then the
claim is proved. Thus we may assume that there are at most two
such edges. Since $P_3\nsubseteq G^g$ and
$P_m\nsubseteq G^b$, there is at least one such an  edge.
Therefore, it is sufficient to consider the following
cases.\\

\noindent$i)$ W.l.g, $G^r$ contains two edges in
$E(\{v_i,v_{i+1}\},\{x,y,z\})$
\vspace{.25 cm}

\noindent$ii)$ W.l.g, $G^r$ contains exactly one edge in
$E(\{v_i,v_{i+1}\},\{x,y,z\})$.\\\\
\noindent If $(i)$ occurs, we may assume that there are exactly
one edge from each of $v_i$ and $v_{i+1}$ to $\{x,y,z\}$ in $G^r$,
otherwise we have nothing to prove. Suppose there is no red edge
in $E(\{v_i,v_{i+1}\},\{z\})$. Since $P_3\nsubseteq G^g$ and
$P_m\nsubseteq G^b$, $G^r$ contains at least one edge in
$E(\{z,u_1,u_2\},\{v_i,v_{i+1}\})$.
Whereas $Q$ is maximal, this edge has to be in
$E(\{v_i,v_{i+1}\},\{z\})$, a contradiction.
\vspace{.5 cm}

\noindent If $(ii)$ occurs, we may assume that
$xv_{i}\in E(G^r)$. Since $G^r$ contains no edge in
$E(\{v_i,v_{i+1}\},\{y,z\})$, green and also blue edges in
$E(\{v_i,v_{i+1}\},\{y,z\})$ form a matching.
Thus, clearly there are two red edges in
$E(\{v_j,v_{j+1}\},\{y,z\})$.
 The reminder of the proof is the same to the case
$(i)$.$\hfill\blacksquare$ \vspace{.6 cm}

 Now, let $Q'$ be a maximal path in the
subgraph of $G^r$ induced by the edges in $E(A\backslash V(Q),K)$
with endpoints   $w_1$ and $w_2$ in $A\backslash V(Q)$ and
$K'=K\backslash V(Q')$.

\vspace{.6 cm} \noindent\textbf{Case 1. }$|A\setminus (V(Q)\cup
V(Q'))|=0$.\vspace{.3 cm}

\noindent Using Lemma \ref{independent red edges}, we may assume
that $G^r$ contains a cycle $C=w_1Q'w_2v_{m-1}u_2Qu_1v_1w_1$ on
 $2\lfloor \frac{n}{2}\rfloor$ vertices. If $n$ is even, we are done. Otherwise, since
$|K'|\geq \lceil\frac{n}{2}\rceil-1$, there is one pair of
vertices $(v_i,v_{i+1})$ in $K'$. Since $G^r$
contains at least one edge in
$E(\{u_1,u_2,w_1,w_2\},\{v_i,v_{i+1}\})$,
 we may  suppose that $v_iu_1\in E(G^r)$ and so
 $R'=v_iu_1Qu_2v_{m-1}w_2Q'w_1v_1$ is a red $P_n$.

\vspace{.6 cm} \noindent\textbf{Case 2. }$|A\setminus (V(Q)\cup
V(Q'))|=1$.\vspace{.3 cm}

\noindent Let $A\setminus (V(Q)\cup V(Q'))=\{x\}$. Using Lemma
\ref{independent red edges} we may assume that
$u_1v_1,u_2v_{m-1},w_1v_{m-1}$ and $w_2v_1$ are red edges. Since
$P_3\nsubseteq G^g$, $G^r$ contains at
 least one edge in $E(\{v_1,v_{m-1}\},\{x\})$, say $xv_1$. Thus
$R=xv_1u_1Qu_2v_{m-1}w_1Q'w_2$ is a red $P_{2\lfloor
\frac{n}{2}\rfloor-1}$.
We consider the following subcases.\\

\noindent {\it Subcase 1.} $n$ is even:

\medskip\noindent Since $|K'|\geq
\lceil\frac{n}{2}\rceil$, there is at least one pair of vertices
 $(v_i,v_{i+1})$ in $K'$.
 If $xv_i\ ($or
$xv_{i+1}$) is red, then $v_ixRw_2$ (or $v_{i+1}xRw_2$
) form a red $P_n$. Otherwise, we may assume that $xv_i\in
E(G^b)$ and $xv_{i+1}\in E(G^g)$. Therefore $xv_{m-1}\in E(G^r)$
and $R'=u_1Qu_2v_{m-1}xv_1w_2Q'w_1$ is a red $P_{n-1}$. Whereas
$G^r$ contains at least one edge of
$E(\{v_i,v_{i+1}\},\{u_1,w_1\})$,
we can extend $R'$ to a red $P_n$.\\

\noindent {\it Subcase 2.} $n$ is odd:

\medskip\noindent Since $|K'|\geq
\lceil\frac{n}{2}\rceil=\frac{n+1}{2}$, there are at least two
disjoint pairs of vertices $(v_j,v_{j+1})$ and $(v_k,v_{k+1})$ in
$K'$. Clearly, each of $x$ and $w_2$ in $G^r$ has at least one
neighbor in $B=\{v_j,v_{j+1},v_k,v_{k+1}\}$, say $s_1$ and $s_2$
respectively. If $s_1\neq s_2$, $s_1xRw_2s_2$ is a red $P_n$, else
$s_1xRw_2s_1$ is a red $C_{n-1}$. One can easily check that $G^r$
contains at least one edge of $E(B\backslash
\{s_1\},\{u_1,u_2,w_1\})$, and so adding
this edge to $C_{n-1}$ yields a $P_n\subseteq G^r$.

\vspace{.6 cm} \noindent\textbf{Case 3. }$|A\setminus (V(Q)\cup
V(Q'))|\geq 2$.\vspace{.3 cm}

\noindent Let $x,y\in A\setminus (V(Q)\cup V(Q'))$. We show that
this case is impossible. Since $|K'|\geq
\lceil\frac{n}{2}\rceil+1$ and at most
$\lfloor\frac{n}{2}\rfloor-4$ vertices of $V(P)\setminus \{v_1,v_{m-1}\}$ are covered by $Q$
and $Q'$, by the Pigeonhole principle we have one of the following
cases.\\\\
\noindent$i)$ $K'$ contains four disjoint pairs of vertices
$(v_k,v_{k+1})$, $(v_i,v_{i+1})$, $(v_j,v_{j+1})$ and
$(v_l,v_{l+1})$.
\vspace{0.25 cm}

\noindent$ii)$ $K'$ contains three consecutive vertices $v_k,
v_{k+1},v_{k+2}$.\\\\
If $(i)$ occurs, since $P_3\nsubseteq G^g$ and $P_m\nsubseteq G^b$ there is a red edge between $x$ and any two
pairs of vertices and so w.l.g we may assume that $xv_{k+1}$, $xv_{l+1}$,
$xv_{i+1}\in E(G^r)$. Since $Q$ and $Q'$ are maximal, $G^r$ contains no
edge in $E(\{u_1,u_2,w_1,w_2\},\{v_{k+1},v_{l+1},v_{i+1}\})$. If
there is a red edge in $E(\{u_1,u_2\},\{v_t,v_{t+1}\})$ (resp. in $E(\{w_1,w_2\},\{v_t,v_{t+1}\})$) for some
$t\in \{i,k,l\}$, then the maximality of $Q$ and $Q'$ implies that green and also blue edges in
$E(\{w_1,w_2\},\{v_t,v_{t+1}\})$ (resp. in $E(\{u_1,u_2\},\{v_t,v_{t+1}\})$) form perfect matchings on four vertices.
Now since there is at least one red edge in $E(\{u_1,u_2\},\{v_t,v_{t+1}\})$ (resp. in $E(\{w_1,w_2\},\{v_t,v_{t+1}\})$) for some
$t\in \{i,k,l\}$, w.l.g we may assume that green and also blue edges in
both $E(\{u_1,u_2\},\{v_k,v_{k+1}\})$ and
$E(\{w_1,w_2\},\{v_l,v_{l+1}\})$ form matchings.
Therefore $\{u_1v_{i+1}, u_2v_{i+1}, w_1v_{i+1},
w_2v_{i+1}\}\subseteq E(G^b)$ and consequently $\{u_1v_i, u_2v_{i},
w_1v_i, w_2v_i\}\subseteq E(G^r)$ which is a contradiction.

\vspace{0.6 cm}

\noindent If $(ii)$ occurs,  at least five vertices of
$\{u_1,u_2,w_1,w_2,x,y\}$ are adjacent to some vertices of $\{v_k,
v_{k+1}, v_{k+2}\}$ in $G^r$, since $P_3\nsubseteq G^g$ and $P_m\nsubseteq G^b$. Let $B$ be the set of the vertices
in $\{u_1,u_2,w_1,w_2,x,y\}$ that are adjacent to a vertex in
$\{v_k,v_{k+1},v_{k+2}\}$ by a red edge. Since $P_3\nsubseteq G^g$ and $P_m\nsubseteq G^b$ then every vertex of $B$ has exactly one red neighbor in $\{v_k, v_{k+1}, v_{k+2}\}$.  Now, we have the following subcases.\\

\noindent {\it Subcase 1.} $\{x,y\}\subseteq B$:

\medskip\noindent By the maximality of $Q$ and $Q'$, we may suppose that
the edges $xv_{t},yv_{t},w_1v_{t'},w_2v_{t'}$ are red for some $t,t'\in\{k,k+1,k+2\}$, $t<t'$ and $u_1v_r\in E(G^r)$ where $r\neq t,t'$.
If $t,t'\in\{k,k+1\}$ (resp. $t,t'\in\{k+1,k+2\}$) then green and also blue edges in
$E(\{x,y\},\{v_{r},v_{t'}\})$ (resp. $E(\{w_1,w_2\},\{v_{r},v_{t}\})$) form matchings
and so there is a red edge in $E(\{u_1,u_2\},\{v_{t},v_{t'}\})$ (resp. $E(\{u_1,u_2\},\{v_{t},v_{t'}\})$), and
this contradicts the maximality of $Q$ and $Q'$. Finally if $t,t'\in\{k,k+2\}$ then green and also blue edges in
$E(\{w_1,w_2\},\{v_{r},v_{t}\})$ form matchings
and so there is a red edge in $E(\{x,y\},\{v_{r},v_{t'}\})$ and again
this contradicts the maximality of $Q$ and $Q'$.\\

\noindent {\it Subcase 2.} $\{x,y\}\cap B=\{x\}$:

\medskip\noindent By a similar argument as in subcase 1, we have a contradiction which completes the proof of the theorem.
$\hfill \blacksquare$
\\\\
{\bf Proof of Theorem \ref{main result}.} It is clear that $5\leq
R(P_3,P_3,P_3)\leq R(P_3,P_3,P_4)$. On the other hand by corollary
\ref{extermal P5,P4}, $R(P_3,P_3,P_4)\leq 5$. Then
$R(P_3,P_3,P_3)=R(P_3,P_3,P_4)=5$.
 Combining Theorems
\ref{R(P3,P4,Pm) and R(P3,P5,Pm)}, \ref{R(P3,P3,Pm)},
\ref{R(p3,Pm,P6),R(P3,Pm,P7)} and \ref{R(P3,Pm,Pn)} give a proof
for Theorem \ref{main result}. $\hfill \blacksquare$

\begin{corollary}
$R(P_3,nK_2,mK_2)=2m+n-1$ for every $m\geq n\geq 3$.
\end{corollary}
\textbf{Proof. }To see $2m+n-1\leq R(P_3,nK_2,mK_2)$, let
$H=K_{n-1}+{\bar K_{2m-1}}$ and ${\bar H}$ be the complement of
$H$ with respect to $K_{2m+n-2}$.
 Clearly coloring $H$ by red and ${\bar H}$ by blue yields a 2-edge coloring of $K_{2m+n-2}$ such that $nK_2\nsubseteq G^r$ and $mK_2\nsubseteq G^b$.
  This means that $2m+n-1\leq R(P_3,nK_2,mK_2)$. Now we prove the upper bound. It is easy to see that $R(P_3,nK_2,mK_2)\leq R(P_3,P_{2n},P_{2m})$
  and by Theorem \ref{main result}, $R(P_3,P_{2n},P_{2m})=2m+n-1$. This observation completes the proof.

$\hfill \blacksquare$


\footnotesize
\begin{thebibliography}{99}

\bibitem{Boundy}
J. A. Bondy, U. S. R. Murty, Graph Theory With Applications, American Elsevier Publishing Co. INC, 1976.

\bibitem{Ramsey number of stripes}
E. J. Cockayne, P. J. Lorimer, The Ramsey number for stripes, {\it J. Austral. Math. Soc.} 19 (Series A) (1975), 252-256.

\bibitem{Electron. J. Combin}
T. Dzido, M. Kubale, K. Piwakowski, On some Ramsey and Tur\'{a}n-type numbers for paths and cycles, {\it Electron. J. Combin.} \#R55 13 (2006).

\bibitem{Discuss. Math. Graph Theory}
T. Dzido, Multicolor Ramsey numbers for paths and cycles, {\it Discuss. Math. Graph Theory} 25 (2005) 57-65.



\bibitem{Path Ramsey numbers in multicolorings}
R. J. Faudree, R. H. Schelp, Path Ramsey numbers in multicolorings, {\it J. Combin. Theory, Ser. B} 19 (1975), 150-160.

\bibitem{Ramsey number for a triple of long even cycles}
A. Figaj, T. Luczak, The Ramsey number for a triple of long even cycles,  {\it J. Combin. Theory, Ser. B} 97 (2007), 584-596.

\bibitem{Ramsey number of paths}
L. Gerencs\'{e}r, A. Gy\'{a}rf\'{a}s, On Ramsey-Type Problems,
{\it Ann. Univ. Sci. Budapest. E\"{o}tv\"{o}s Sect. Math. } 10
(1967), 167-170.

\bibitem{Ramsey number of three-paths}
 A. Gy\'{a}rf\'{a}s, M. Ruszink\'{o}, G. S\'{a}rk\"{o}zy, E. Szemer\'{e}di, Three-color Ramsey numbers for paths, {\it Combinatorica} 27 (1) (2007), 35-69.

\bibitem{tripartite Ramsey number of paths}
A. Gy\'{a}rf\'{a}s, M. Ruszink\'{o}, G. S\'{a}rk\"{o}zy, E.
Szemer\'{e}di, Tripartite Ramsey numbers for paths, {\it J. Graph
Theory} 55 (2007), 164-170.

\bibitem{Omidi Raeisi}
G.R. Omidi, G. Raeisi, On multicolor Ramsey number of paths versus
cycles, Electron. J. Combin. \#P24 18 (2011).

\bibitem{survey}
 S. P. Radziszowski, Small Ramsey numbers, {\it Electron. J. Combin.} 1 (1994), Dynamic Surveys, DS1.12 (August 4, 2009).

\bibitem{Ramsey}
F. P. Ramsey, On a problem of formal logic, {\it Proc. London Math. Soc. 2nd Ser.} 30 (1930), 264-286.

\bibitem{schelp conj}
R. H. Schelp, Some Ramsey-Turan type problems and related questions,
manuscript.

\bibitem{schelp}
R. H. Schelp, A minimum degree condition on a Ramsey graph which arrows a
path, manuscript.

\bibitem{Euro. J. Combin}
Z. Shao, X. Xu, X. Shi, L. Pan, Some three-color Ramsey numbers,
$R(P_4,P_5,C_k)$ and $R(P_4,P_6,C_k)$, {\it Europ. J. Combin.} 30
(2009), 396-403.

\end{thebibliography}
\end{document}